\documentclass[11pt,a4paper]{article}

\usepackage[utf8]{inputenc}
\usepackage[T1]{fontenc}
\usepackage{textcomp}
\usepackage{amsmath,amssymb,amsthm}
\usepackage{hyperref}
\usepackage{geometry}
\usepackage{microtype}
\microtypesetup{expansion=false}
\usepackage{verbatim}
\usepackage{enumitem}
\usepackage{xcolor}
\definecolor{codegray}{gray}{0.95}

\theoremstyle{plain}
\newtheorem{theorem}{Theorem}[section]

\theoremstyle{definition}

\newcommand{\Z}{\mathbb{Z}}
\newcommand{\N}{\mathbb{N}}

\newcommand{\Heis}{\operatorname{Heis}}

\newcommand{\mathlib}{\textsf{mathlib}}
\newcommand{\code}[1]{\texttt{#1}}

\title{\bfseries Classifying the Groups of Order $p^3$ in Lean}
\author{
  Li Xiang\thanks{The formalisation and this manuscript were prepared with assistance from AI language models. The author reviewed and validated all AI-generated content.}
}
\date{June 2026}

\begin{document}
\maketitle

\begin{abstract}
This note discusses our formalisation in Lean 4 of the classification of groups of order $p^3$ for a prime number $p$, using mathlib4. We present the five isomorphism classes and give a detailed account of the formalisation, with particular emphasis on the non-abelian case, which requiring the most substantial formal development. For odd~$p$, the non-abelian groups are the Heisenberg group $\Heis(\Z/p\Z)$ and the semidirect product $\Z/p^2\Z\rtimes\Z/p\Z$; for $p=2$, they are $D_4$ and $Q_8$. We describe the construction of these concrete groups, the structural lemmas about centers, commutators, and exponents, and the explicit isomorphism constructions that classify an arbitrary non-abelian $p^3$-group.
\end{abstract}

\section{Introduction}

The computer-assisted formalisation of mathematics sees proofs transcribed into computer code and rigorously checked by software. The formal verification of the Feit--Thompson Odd Order Theorem~\cite{gonthier2013} using the Coq proof assistant drew great attention from the mathematical community. Lean has rapidly become one of the most actively used theorem provers, with an extensive mathematical library \mathlib{}~\cite{mathlib2020}. An influential theme in finite group theory is classification. While the Classification of Finite Simple Groups is far out of reach, many smaller classification results are accessible and serve as stepping stones.

All groups of prime order are cyclic (already formalised in Lean). Classifications of groups of order $p^2$ and $pq$ (for primes $p,q$) have been formalised by Harper and Wu~\cite{harper2025pq}. Our work extends this to groups of order $p^3$.

\textbf{Formalisation.} The classification of groups of order $p^3$ for an arbitrary prime $p$.

This work contains no new mathematical results. The classification of groups
of order $p^3$ is a classical theorem found in standard group theory
textbooks. Our contribution lies in the formal verification and
organisation of this classification in Lean~4, making it the first
complete formalisation of this result in any proof assistant.

There are five isomorphism classes. Three are abelian:
\[
\Z/p^3\Z,\qquad \Z/p^2\Z \times \Z/p\Z,\qquad (\Z/p\Z)^3.
\]
Two are non-abelian. For odd~$p$:
\begin{itemize}[nosep]
  \item $\Heis(\Z/p\Z)$, of exponent $p$;
  \item $\Z/p^2\Z \rtimes \Z/p\Z$, of exponent $p^2$.
\end{itemize}
For $p=2$, these become the dihedral group $D_4$ and the quaternion group $Q_8$ (both of order~$8$).

The formalisation consists of approximately 3000 lines of Lean code in six files: \code{Defs.lean}, \code{Structural.lean}, \code{AbelianCase.lean}, \code{NonAbelianCase.lean}, \code{Classification.lean}, and the top-level \code{P3Group.lean}. Of these, \code{NonAbelianCase.lean} (over 2000 lines) contains the most intricate arguments.

The full Lean code is available at \url{https://github.com/lixiang90/p3group}.

\subsection{Contributions}

\begin{itemize}[nosep]
  \item First complete formalisation of the $p^3$ classification in any proof assistant (\textasciitilde 3000 lines of Lean~4, verified by the Lean kernel).
  \item Construction of the five canonical groups of order $p^3$ in Lean, including the Heisenberg group and the semidirect product $\Z/p^2\Z\rtimes\Z/p\Z$.
  \item A library of structural lemmas for non-abelian $p^3$-groups (center, commutator, nilpotency class, exponent dichotomy).
  \item Explicit isomorphism constructions via the structure theorem, the Hall--Petrescu formula, and case analysis on the exponent and $p=2$.
\end{itemize}

\section{Underlying mathematics}

\subsection{Abelian groups of order \texorpdfstring{$p^3$}{p\^{}3}}

Let $G$ be a finite abelian group of order $p^3$. By the structure theorem for finite abelian groups, $G$ is a direct product of cyclic groups whose orders are prime powers multiplying to $p^3$. The partitions of~$3$ correspond to the three abelian types:
\[
\Z/p^3\Z,\qquad \Z/p^2\Z \times \Z/p\Z,\qquad (\Z/p\Z)^3.
\]
These are pairwise non-isomorphic: the first is cyclic, the second is non-cyclic with exponent $p^2$, and the third has exponent $p$.

\subsection{Non-abelian groups of order \texorpdfstring{$p^3$}{p\^{}3}}

\begin{theorem}[Structure of non-abelian $p^3$-groups]\label{thm:structure}
Let $G$ be a non-abelian group of order $p^3$ with $p$ prime.
\begin{enumerate}[label=(\alph*),nosep]
  \item $|Z(G)| = p$ and $G/Z(G) \cong (\Z/p\Z)^2$.
  \item The commutator subgroup $[G,G]$ equals $Z(G)$.
  \item $G$ has nilpotency class exactly~$2$.
  \item $\exp(G)$ is either $p$ or $p^2$.
\end{enumerate}
\end{theorem}

\noindent \textbf{Exponent $p$ (odd $p$).} When $\exp(G)=p$ and $p \neq 2$, $G \cong \Heis(\Z/p\Z)$. The Heisenberg group is the set of triples $(a,b,c) \in (\Z/p\Z)^3$ with multiplication
\[
(a,b,c)\cdot(a',b',c') = (a+a',\, b+b',\, c+c'+ab').
\]

\noindent \textbf{Exponent $p^2$ (odd $p$).} When $\exp(G)=p^2$ and $p \neq 2$, $G \cong \Z/p^2\Z \rtimes \Z/p\Z$, presented as
\[
\langle\, a,b \mid a^{p^2}=b^p=1,\; b^{-1}ab = a^{1+p} \,\rangle.
\]

\noindent \textbf{The case $p=2$.} The non-abelian groups of order $8$ are $D_4$ and $Q_8$, distinguished by the number of involutions: $D_4$ has five elements of order $\leq 2$, $Q_8$ has only one.

\section{Main statements in Lean}

We briefly describe the Lean formalisation of the key definitions and theorems.
Group structures are formalised by the \code{Group} type class; \code{G ->* H}
and \code{G $\simeq$* H} are the types of homomorphisms and isomorphisms;
\code{Nat.card G} gives the cardinality of a type.

\subsection{The five standard models}

We define five concrete groups of order $p^3$, each serving as a canonical
representative of one isomorphism class.  The three abelian models are built
from \code{ZMod} (the multiplicative cyclic group $\Z/n\Z$):
\begin{align*}
  C_{p^3} &= \Z/p^3\Z, &
  C_{p^2}\times C_p &= \Z/p^2\Z \times \Z/p\Z, &
  (C_p)^3 &= (\Z/p\Z)^3.
\end{align*}
These appear in Lean as \code{CyclicP3}, \code{AbelianP2P}, and
\code{ElementaryP3}.

For odd~$p$, the two non-abelian models are custom structures.
The \textbf{Heisenberg group} $\Heis(\Z/p\Z)$ consists of triples
$(a,b,c)\in(\Z/p\Z)^3$ with multiplication
\[
(a,b,c)\cdot(a',b',c')=(a+a',\;b+b',\;c+c'+ab').
\]
The \textbf{semidirect product} $\Z/p^2\Z\rtimes\Z/p\Z$ has elements
$(a,b)$ with $a\in\Z/p^2\Z$, $b\in\Z/p\Z$ and multiplication
\[
(a,b)\cdot(a',b')=\bigl(a+a'+b_{\mathrm{val}}\!\cdot p\cdot a',\;b+b'\bigr),
\]
where $b_{\mathrm{val}}$ is the natural-number lift of $b$.
These are defined in Lean as \code{HeisenbergGroup p} and
\code{SemidirectP2P p}.

For $p=2$, the two non-abelian groups of order~$8$ are the dihedral group
$D_4$ and the quaternion group $Q_8$, both already available in \mathlib{}
as \code{DihedralGroup 4} and \code{QuaternionGroup 2}.

\subsection{The classification predicate}

The property ``$G$ is isomorphic to one of the five standard groups of
order $p^3$'' is captured by the predicate
\[
\operatorname{IsP3Group}_p(G)\;\longleftrightarrow\;
\begin{array}{l}
G\cong C_{p^3}\;\lor\;
G\cong C_{p^2}\times C_p\;\lor\;
G\cong (C_p)^3\;\lor\\
(p\neq2\;\land\;G\cong\Heis(\Z/p\Z))\;\lor\;
(p\neq2\;\land\;G\cong\Z/p^2\Z\rtimes\Z/p\Z)\;\lor\\
(p=2\;\land\;G\cong D_4)\;\lor\;
(p=2\;\land\;G\cong Q_8).
\end{array}
\]
In Lean this is the inductive predicate \code{IsP3Group p G}; see
\code{P3Group/Classification.lean} for its definition.

The central theorem of the formalisation is then:

\begin{theorem}[Classification of groups of order $p^3$]
Let $p$ be prime and $G$ a finite group.  If $|G|=p^3$, then
$\operatorname{IsP3Group}_p(G)$.
\end{theorem}

In Lean this reads
\code{theorem classification (G : Type*) [Group G]}
\code{[Fintype G] (hcard : Nat.card G = p \^{} 3) :}
\code{IsP3Group p G}.

\subsection{Pairwise non-isomorphism}

To complete the classification, we prove that the five standard models are
mutually non-isomorphic:

\begin{itemize}[nosep]
  \item $C_{p^3}\not\cong C_{p^2}\times C_p$ and
    $C_{p^3}\not\cong(C_p)^3$ because the former is cyclic while the
    latter two are not;
  \item $C_{p^2}\times C_p\not\cong(C_p)^3$ because their exponents differ
    ($p^2$ versus $p$);
  \item $\Heis(\Z/p\Z)\not\cong\Z/p^2\Z\rtimes\Z/p\Z$ for odd~$p$, again
    distinguished by exponent ($p$ versus $p^2$);
  \item $D_4\not\cong Q_8$, distinguished by the number of involutions
    ($D_4$ has five, $Q_8$ has one).
\end{itemize}

The corresponding Lean theorems are
\code{abelian\_types\_distinct},
\\ \code{heisenberg\_not\_iso\_semidirect}, and
\code{dihedral4\_not\_iso\_quaternion8}.

\subsection{Structural lemmas for the non-abelian case}

The non-abelian classification depends on the following structural facts,
all proved in \\ \code{P3Group/Structural.lean}.
For a non-abelian group $G$ with $|G|=p^3$:
\begin{enumerate}[label=(\alph*),nosep]
  \item $|Z(G)|=p$ and $G/Z(G)\cong(\Z/p\Z)^2$;
  \item $[G,G]=Z(G)$;
  \item the nilpotency class of $G$ is exactly~$2$;
  \item $\exp(G)\in\{p,p^2\}$.
\end{enumerate}
The corresponding Lean statements are
\code{center\_card\_eq\_p\_of\_nonabelian},
\\ \code{quotient\_center\_iso\_p2},
\code{commutator\_eq\_center},
\\ \code{nilpotencyClass\_eq\_two}, and
\code{exponent\_of\_nonabelian\_p3}.

\section{Implementation details of the abelian case}

This section describes the proof in \code{AbelianCase.lean}. The file is concise (\textasciitilde 240 lines) yet contains a complete argument using the structure theorem for finite abelian groups.

\subsection{Concrete abelian groups and their invariants}

The three abelian types are defined as abbreviations:

\begin{verbatim}
abbrev CyclicP3 := ZMod (p ^ 3)
abbrev AbelianP2P := ZMod (p ^ 2) x ZMod p
abbrev ElementaryP3 := ZMod p x ZMod p x ZMod p
\end{verbatim}

We verify each has cardinality $p^3$ (\code{card\_cyclicP3}, \code{card\_abelianP2P},
\code{card\_elementaryP3}). To prove the three are pairwise non-isomorphic,
we compute their distinguishing invariants:

\begin{verbatim}
theorem abelianP2P_not_cyclic :
    ~ IsCyclic (Multiplicative (ZMod (p ^ 2)) x
                Multiplicative (ZMod p))
\end{verbatim}

This uses \mathlib{}'s lemma \code{coprime\_card\_of\_isCyclic\_prod}: a product of two cyclic groups is cyclic only if their orders are coprime. Here $p^2$ and $p$ are not coprime.

\begin{verbatim}
theorem elementaryP3_exponent :
    exponent (Multiplicative (ZMod p) x Multiplicative (ZMod p) x
              Multiplicative (ZMod p)) = p

theorem abelianP2P_exponent :
    exponent (Multiplicative (ZMod (p ^ 2)) x Multiplicative (ZMod p)) = p ^ 2
\end{verbatim}

These use \code{exponent\_prod} and \code{exponent\_multiplicative} from \mathlib{},
reducing to computing the lcm of the exponents of each factor: $\operatorname{lcm}(p,p,p)=p$
and $\operatorname{lcm}(p^2,p)=p^2$.

\subsection{The structure theorem and the main proof}

The central result is \code{abelian\_p3\_classification}:

\begin{verbatim}
theorem abelian_p3_classification (G : Type*) [CommGroup G] [Fintype G]
    (hcard : Nat.card G = p ^ 3) :
    Nonempty (G ~=* Multiplicative (CyclicP3 p)) \/
    Nonempty (G ~=* (Multiplicative (ZMod (p ^ 2)) x
                     Multiplicative (ZMod p))) \/
    Nonempty (G ~=* (Multiplicative (ZMod p) x
                     Multiplicative (ZMod p) x
                     Multiplicative (ZMod p)))
\end{verbatim}

The proof proceeds in two main cases.

\subsubsection{Case 1: \texorpdfstring{$G$}{G} is cyclic}

If $G$ is cyclic, then $G \cong \Z/|G|\Z = \Z/p^3\Z$. In Lean we use \code{zmodCyclicMulEquiv},
which gives an isomorphism $G \cong \Z/|G|\Z$ in multiplicative form
(\code{Multiplicative (ZMod (Nat.card G))}). Since $\operatorname{Nat.card}(G) = p^3$,
this is exactly the first isomorphism type.

\subsubsection{Case 2: \texorpdfstring{$G$}{G} is not cyclic}

If $G$ is not cyclic, we invoke \mathlib{}'s structure theorem for finite abelian groups:

\begin{verbatim}
obtain <i, inst, n, hn_gt, <e>> :=
    CommGroup.equiv_prod_multiplicative_zmod_of_finite G
\end{verbatim}

This yields an index type $\iota$, a sequence $n_i$ of natural numbers with $n_i > 1$,
and an isomorphism $e : G \to \prod_i \Z/n_i\Z$ (multiplicative).

\textbf{Step A: Each $n_i$ is a power of $p$.}
Since $\prod_i n_i = p^3$ (by cardinality), each $n_i \mid p^3$.
Because $p$ is prime, $n_i = p^{k_i}$ for some $1 \leq k_i \leq 3$.

\textbf{Step B: $2 \leq |\iota| \leq 3$.}
If $|\iota| = 0$, the product is empty (order $1$, contradicting $p^3 > 1$).
If $|\iota| = 1$, then $G$ would be cyclic, contradicting the assumption.
Hence $|\iota| \geq 2$.
On the other hand, $p^{|\iota|} \leq \prod_i n_i = p^3$ (since each $n_i \geq p$),
so $|\iota| \leq 3$. Thus $|\iota| \in \{2, 3\}$.

\textbf{Step C: $|\iota| = 2$.}
Reindex $\iota$ to \code{Fin 2} via \code{Fintype.equivFinOfCardEq}.
Let the two factors be $p^{k_0}$ and $p^{k_1}$ with $k_0, k_1 \geq 1$ and $k_0 \leq k_1 \leq 3$.
Since $k_0 + k_1 = 3$, the only possibilities are $(k_0,k_1) = (1,2)$.
Thus $G \cong \Z/p^2\Z \times \Z/p\Z$. If the order is $(p, p^2)$ we compose
with \code{MulEquiv.prodComm} to swap to $(p^2, p)$.

\textbf{Step D: $|\iota| = 3$.}
Reindex to \code{Fin 3}. Each $k_i \geq 1$ and $k_0+k_1+k_2 = 3$, so
$k_0 = k_1 = k_2 = 1$. Thus all factors are $p$, giving $G \cong (\Z/p\Z)^3$.
We use an explicit \code{MulEquiv} from $\prod_{j:\mathrm{Fin}\,3} \Z/p\Z$ to
$(\Z/p\Z) \times (\Z/p\Z) \times (\Z/p\Z)$ defined by $f \mapsto (f(0), f(1), f(2))$.

\subsection{Helper lemmas for the product-to-tuple conversion}

Two private definitions bridge the gap between $\prod_{i:\iota} \Z/n_i\Z$ (an indexed product)
and the concrete product types used in the classification:

\begin{verbatim}
private def mulEquivPiReindex {i i' : Type*} [Fintype i] [Fintype i']
    [DecidableEq i'] (M : i -> Type*) [forall i, Mul (M i)] (e : i = i') :
    (forall i : i, M i) ~=* (forall j : i', M (e.symm j))

private def mulEquivPiFinTwo (M : Fin 2 -> Type*) [forall i, Mul (M i)] :
    (forall i : Fin 2, M i) ~=* M 0 x M 1
\end{verbatim}

The first reindexes the product along an equivalence of index sets; the second converts
a 2-element indexed product to a binary product type---exactly \code{piFinTwoEquiv M}
extended to respect multiplication.

\section{Implementation details of the non-abelian case}

This section describes the proof architecture of \code{NonAbelianCase.lean} in detail. The overall strategy is:

\begin{enumerate}[nosep]
  \item Construct the two concrete non-abelian groups (Heisenberg and semidirect product) and establish their basic properties (order, non-abelianness, exponent).
  \item Prove the exponent dichotomy: a non-abelian $p^3$-group has exponent $p$ or $p^2$.
  \item For exponent $p$ and $p \neq 2$, build an isomorphism $G \cong \Heis(\Z/p\Z)$.
  \item For exponent $p^2$ and $p \neq 2$, build an isomorphism $G \cong \Z/p^2\Z\rtimes\Z/p\Z$.
  \item For $p=2$, build an isomorphism to $D_4$ or $Q_8$.
\end{enumerate}

\subsection{The Heisenberg group}

The Heisenberg group $\Heis(\Z/p\Z)$ is a structure with three coordinates in \code{ZMod p}. Its group instance is given by:
\[
(x \cdot y).a = x.a + y.a,\qquad
(x \cdot y).b = x.b + y.b,\qquad
(x \cdot y).c = x.c + y.c + x.a \cdot y.b.
\]
The inverse is $(x^{-1}).a = -x.a$, $(x^{-1}).b = -x.b$, $(x^{-1}).c = -x.c + x.a \cdot x.b$.

\subsubsection{Power formula}

A key technical lemma gives an explicit formula for powers in the Heisenberg group, proved by induction on $n$:
\[
(x^n).c = n \cdot x.c + \binom{n}{2} \cdot (x.a \cdot x.b) \qquad\text{(in $\Z/p\Z$)}.
\]

The full Lean statement reads:

\begin{verbatim}
private theorem heisenberg_pow_aux (p : Nat) (x : HeisenbergGroup p) (n : Nat) :
    (x ^ n).a = (n : ZMod p) * x.a /\
    (x ^ n).b = (n : ZMod p) * x.b /\
    (x ^ n).c = (n : ZMod p) * x.c +
      (Nat.choose n 2 : ZMod p) * (x.a * x.b)
\end{verbatim}

When $n = p$ and $p \neq 2$, we have $(p : \Z/p\Z) = 0$ and $\binom{p}{2} \equiv 0 \pmod p$, so $x^p = 1$ for all $x$. This gives the upper bound $\exp(G) \mid p$. To see $\exp(G) = p$, we exhibit the element $(1,0,0)$ of order exactly $p$ using the same power formula. Together, $\exp(G) = p$.

\subsection{The semidirect product \texorpdfstring{$\Z/p^2\Z \rtimes \Z/p\Z$}{Z/p2Z x Z/pZ}}

This group is a structure with coordinates $a \in \Z/p^2\Z$, $b \in \Z/p\Z$. The multiplication is:
\[
(a,b)\cdot(a',b') = \bigl(a + a' + b_{\text{val}}\cdot p \cdot a',\; b + b'\bigr),
\]
where $b_{\text{val}}$ denotes the natural number lift of $b$.

The group axioms require a number of ring arithmetic lemmas in $\Z/p^2\Z$, most notably:
\begin{itemize}[nosep]
  \item $p \cdot p = 0$ in $\Z/p^2\Z$ (\code{pp\_eq\_zero}).
  \item $((k_1+k_2).\mathrm{val} : \Z/p^2\Z)\cdot p = (k_1.\mathrm{val})\cdot p + (k_2.\mathrm{val})\cdot p$ (\code{val\_mul\_p\_add}).
  \item $((-k).\mathrm{val} : \Z/p^2\Z)\cdot p = -(k.\mathrm{val})\cdot p$ (\code{val\_neg\_mul\_p}).
\end{itemize}

\subsubsection{Exponent}

For the exponent proof, we need formulas for powers in the semidirect product. One key step:
\begin{itemize}[nosep]
  \item If $y.b = 0$, then $(y^k).a = k \cdot y.a$ (the $b$ component does not contribute).
  \item When $x.b = 0$, $(x^p).b = 0$. Then $(x^p)^p = (x^p.a)^p$ in the first component, which vanishes using $p\cdot p = 0$.
\end{itemize}
This proves $x^{p^2} = 1$ for all $x$, so $\exp(G) \mid p^2$. The element $(1,0)$ has order $p^2$, giving equality.

\subsection{The \texorpdfstring{$p=2$}{p=2} models: \texorpdfstring{$D_4$}{D4} and \texorpdfstring{$Q_8$}{Q8}}

For the $p=2$ case we reuse \mathlib{}'s existing constructions. We verify:
\begin{itemize}[nosep]
  \item $\operatorname{Nat.card}(D_4) = 8$, $D_4$ is non-abelian.
  \item $\operatorname{Nat.card}(Q_8) = 8$, $Q_8$ is non-abelian.
  \item $D_4 \not\cong Q_8$. Every element of $D_4$ of the form $\mathrm{sr}(i)$ squares to~$1$; there are four such. In $Q_8$, the only elements squaring to~$1$ are $1$ and $a^2$. An isomorphism would inject these four distinct elements into a set of size~$2$, impossible.
\end{itemize}

\subsection{Structural lemmas}

\subsubsection{Center and quotient}

Since $G$ is a $p$-group and non-abelian, $Z(G) \neq 1$ and $Z(G) \neq G$. By Lagrange's theorem, $|Z(G)| \mid |G| = p^3$, so $|Z(G)| \in \{p, p^2\}$. If $|Z(G)| = p^2$, then $|G/Z(G)| = p$, making the quotient cyclic, which would force $G$ to be abelian---contradiction. Hence $|Z(G)| = p$. Then $|G/Z(G)| = p^2$ by Lagrange's theorem.

To prove $G/Z(G) \cong (\Z/p\Z)^2$, we use the structure theorem for finite abelian groups on the quotient (which is abelian of order $p^2$, since $p$-groups of order $p^2$ are abelian). If the quotient were cyclic, $G$ would be abelian (a group with cyclic central quotient is abelian). Hence the quotient is the non-cyclic group of order $p^2$, i.e.\ $(\Z/p\Z)^2$.

\subsubsection{Commutator and nilpotency class}

The quotient $G/Z(G)$ being abelian means all commutators lie in $Z(G)$, so $[G,G] \subseteq Z(G)$. Since $G$ is non-abelian, $[G,G] \neq 1$. Both subgroups have order $p$, so $[G,G] = Z(G)$. This directly implies nilpotency class $2$ (class $1$ would mean abelian).

\subsubsection{Exponent dichotomy}

The exponent of $G$ divides $|G| = p^3$, so $\exp(G) = p^k$ for some $0 \leq k \leq 3$. It cannot be $1$ (then $G$ is trivial, but $|G|=p^3>1$) and cannot be $p^3$ (an element of order $p^3$ would make $G$ cyclic, hence abelian). Thus $k \in \{1,2\}$.

\subsection{Exponent \texorpdfstring{$p$}{p} case: isomorphism to \texorpdfstring{$\Heis(\Z/p\Z)$}{Heis(Z/pZ)}}

This is the heart of lemma \code{heisenberg\_of\_exponent\_p}. Given $\exp(G) = p$ and $p \neq 2$:

\subsubsection{Finding the generators}

Since $G$ is non-abelian, pick $x, y \in G$ with $x y \neq y x$. Set
\[
z := x^{-1} y^{-1} x y,
\]
the commutator of $x^{-1}$ and $y^{-1}$. Then $z \neq 1$ (otherwise $x$ and $y$ commute). The key structural lemma \code{commutator\_mem\_center\_of\_p3} (Section~4.9) tells us $z \in Z(G)$. Moreover, from $x y = y x \cdot z$ we obtain the fundamental commutation relation:
\[
x \cdot y = y \cdot x \cdot z \qquad\text{(equation $\star$)}.
\]

Because $\exp(G) = p$, every non-identity element has order exactly $p$.
Thus
$\operatorname{orderOf}(x)=\operatorname{orderOf}(y)
=\operatorname{orderOf}(z)=p$.

\subsubsection{An algebraic identity for commuting elements}

The purely group-theoretic lemma \code{heisenberg\_mul\_identity} is the key computational engine. It states that if $z \in Z(G)$ and $x y = y x z$, then for all natural numbers $a_1,a_2,b_1,b_2,c_1,c_2$:
\[
y^{b_1} x^{a_1} z^{c_1} \cdot y^{b_2} x^{a_2} z^{c_2}
= y^{b_1+b_2} x^{a_1+a_2} z^{a_1 b_2 + c_1 + c_2}.
\]

The proof is a systematic rearrangement of factors using the centrality of $z$ and the commutation relation $x^{a_1} y^{b_2} = y^{b_2} x^{a_1} z^{a_1 b_2}$ (which is proved first in lemma \code{pow\_mul\_pow\_comm}). The exponents match exactly the Heisenberg group multiplication law. Its formal statement is:

\begin{verbatim}
private lemma heisenberg_mul_identity {G : Type*} [Group G] {x y z : G}
    (hcent : z \in Subgroup.center G)
    (hrel : x * y = y * x * z)
    (a1 a2 b1 b2 c1 c2 : Nat) :
    y ^ b1 * x ^ a1 * z ^ c1 * (y ^ b2 * x ^ a2 * z ^ c2) =
    y ^ (b1 + b2) * x ^ (a1 + a2) * z ^ (a1 * b2 + c1 + c2)
\end{verbatim}

\subsubsection{The homomorphism}

Define $f : \Heis(\Z/p\Z) \to G$ by
\[
f(a,b,c) = y^{b.\mathrm{val}} \cdot x^{a.\mathrm{val}} \cdot z^{c.\mathrm{val}}.
\]
In Lean this reads:

\begin{verbatim}
let fFun : HeisenbergGroup p -> G := fun <a, b, c> =>
  y ^ b.val * x ^ a.val * z ^ c.val
let f := MonoidHom.mk' fFun hmul
\end{verbatim}
where the multiplicativity proof \code{hmul} uses
\code{heisenberg\_mul\_identity} together with
\code{pow\_zmod\_add} and \code{pow\_zmod\_mul},
to verify $f(g_1 g_2) = f(g_1) f(g_2)$.

\subsubsection{Injectivity}

To prove $\ker f$ is trivial, suppose $f(a,b,c) = 1$, i.e.\ $y^{b.\mathrm{val}} x^{a.\mathrm{val}} z^{c.\mathrm{val}} = 1$. The product is in $Z(G)$ only if $a = 0$ (using the commutation relation to commute past $y$), then $b = 0$ (commuting past $x$), and finally $c = 0$ (since $z^{c.\mathrm{val}} = 1$ implies $c = 0$ in $\Z/p\Z$ when $\operatorname{orderOf}(z) = p$). Hence $f$ is injective.

\subsubsection{Bijectivity}

Both $G$ and $\Heis(\Z/p\Z)$ have cardinality $p^3$. An injective map between two finite sets of equal size is bijective. Thus $f$ is an isomorphism.

\subsection{Exponent \texorpdfstring{$p^2$}{p2} case: isomorphism to \texorpdfstring{$\Z/p^2\Z \rtimes \Z/p\Z$}{Z/p2Z x Z/pZ}}

This is lemma \code{semidirectP2P\_of\_exponent\_p2}. Given $\exp(G) = p^2$ and $p \neq 2$:

\subsubsection{Finding \texorpdfstring{$x$}{x} of order \texorpdfstring{$p^2$}{p\^{}2}}

If no element has order $p^2$, then every element has order dividing $p$.
Then $\exp(G) \mid p$, contradicting $\exp(G) = p^2$. So there exists $x$
with $\operatorname{orderOf}(x) = p^2$.

\subsubsection{The subgroup \texorpdfstring{$\langle x\rangle$}{<x>} is normal}

$|\langle x\rangle| = p^2$, so by Lagrange its index is $p$.
A classic lemma (in \mathlib{} as
\\ \code{Subgroup.normal\_of\_index\_eq\_minFac\_card})
says: if a subgroup's index equals the smallest prime factor of $|G|$,
the subgroup is normal. Here both equal $p$, so
$\langle x\rangle \trianglelefteq G$.

\subsubsection{Finding \texorpdfstring{$y$}{y} of order \texorpdfstring{$p$}{p} outside \texorpdfstring{$\langle x\rangle$}{<x>}}

Since $|G| = p^3 > p^2 = |\langle x\rangle|$, there exists $w \notin \langle x\rangle$. If $\operatorname{orderOf}(w) \neq p$, then $\operatorname{orderOf}(w) = p^2$ (since order divides exponent $p^2$ and cannot be $1$). In that case $w^p \in \langle x\rangle$ (because in the quotient $G/\langle x\rangle$, which has order $p$, every element has order dividing $p$). Writing $w^p = x^t$, we set $y = x^{-k} w$ where $t = pk$. One checks $y \notin \langle x\rangle$ and (using the Hall--Petrescu formula, Section~4.7) $\operatorname{orderOf}(y) = p$.

\subsubsection{Normalizing the conjugation action}

Because $\langle x\rangle$ is normal, $y x y^{-1} \in \langle x\rangle$.
Write $y x y^{-1} = x^k$ for some $k \in \Z$. Reducing $k$ modulo $p^2$,
we obtain $m \in \N$ with $0 \leq m < p^2$ such that $y x y^{-1} = x^m$.

\textbf{Step 1: $\gcd(m, p^2) = 1$.}
Conjugation preserves order,
so $\operatorname{orderOf}(x^m) = p^2$.
But $\operatorname{orderOf}(x^m) = p^2 / \gcd(p^2, m)$.
Thus $\gcd(p^2, m) = 1$, which implies $\gcd(m, p) = 1$.

\textbf{Step 2: $m \equiv 1 \pmod p$.}
Consider $z = y x y^{-1} x^{-1} = x^m x^{-1} = x^{m-1}$.
By \\ \code{commutator\_mem\_center\_of\_p3}, $z \in Z(G)$.
Since $|Z(G)| = p$, we have $z^p = 1$.
From $x^{(m-1)p} = 1$ and
$\operatorname{orderOf}(x) = p^2$,
we obtain $p^2 \mid (m-1)p$, hence $p \mid m-1$,
i.e.\ $m \equiv 1 \pmod p$.

\textbf{Step 3: Write $m = 1 + ap$, find $r$ with $ra \equiv 1 \pmod p$.}  Since $m \equiv 1 \pmod p$, write $m = 1 + ap$ with $0 \leq a < p$. If $a = 0$, then $m = 1$, so $y x y^{-1} = x$, meaning $x$ and $y$ commute. A centralizer argument shows this would force $y \in \langle x\rangle$ (since $C_G(x)$ has order at least $p^2$, and if it equaled $p^3$ then $x \in Z(G)$, contradicting $|Z(G)| = p$). Hence $a \neq 0$, so $\gcd(a, p) = 1$ and there exists $r$ with $ra \equiv 1 \pmod p$.

\textbf{Step 4: Replace $y$ by $y^r$.}  By iterating the conjugation relation, one obtains $y^r x (y^r)^{-1} = x^{m^r}$. The number-theoretic lemma \code{one\_add\_mul\_p\_pow\_inv} shows $(1+ap)^r \equiv 1 + p \pmod{p^2}$, so $y' := y^r$ satisfies:
\[
y' x (y')^{-1} = x^{1+p}.
\]
Moreover, $y' \notin \langle x\rangle$ and $\operatorname{orderOf}(y') = p$ (since $\gcd(r,p) = 1$). This entire normalization step is encapsulated in the lemma \code{normalize\_conjugation\_to\_one\_add\_p}:

\begin{verbatim}
private lemma normalize_conjugation_to_one_add_p
    {G : Type*} [Group G] [Fintype G]
    {p : Nat} [Fact p.Prime]
    {x y : G}
    (hx : orderOf x = p ^ 2)
    (hy_ord : orderOf y = p)
    (hy_not_mem : y \notin zpowers x)
    (hconj_mem : y * x * y^(-1) \in zpowers x)
    (hcard : Nat.card G = p ^ 3)
    (hnonab : ~ forall a b : G, a * b = b * a) :
    exists y' : G,
      y' \notin zpowers x /\
      orderOf y' = p /\
      y' * x * y'^(-1) = x ^ (1 + p)
\end{verbatim}

\subsubsection{The homomorphism}

Define $f : (\Z/p^2\Z \rtimes \Z/p\Z) \to G$ by
\[
f(a,b) = x^{a.\mathrm{val}} \cdot (y')^{b.\mathrm{val}}.
\]
To verify $f$ is multiplicative, the key identity is:
\[
(y')^{b_1.\mathrm{val}} \cdot x^{a_2.\mathrm{val}} \cdot (y')^{-b_1.\mathrm{val}}
= x^{a_2.\mathrm{val} \cdot (1+p)^{b_1.\mathrm{val}}}.
\]
Using the binomial congruence $(1+p)^n \equiv 1 + np \pmod{p^2}$ (lemma \code{one\_plus\_p\_pow\_mod\_p\_sq}) and the fact that $\operatorname{orderOf}(x) = p^2$, we reduce the exponent modulo $p^2$ to match the semidirect product law.

\subsubsection{Injectivity and bijectivity}

As before, injectivity uses the kernel argument. If $x^{a.\mathrm{val}} (y')^{b.\mathrm{val}} = 1$, then $(y')^{b.\mathrm{val}} \in \langle x\rangle$. Since $y' \notin \langle x\rangle$ and $G/\langle x\rangle$ has order $p$, this forces $b = 0$. Then $x^{a.\mathrm{val}} = 1$ implies $a = 0$ (since $\operatorname{orderOf}(x) = p^2$). Cardinalities match, so $f$ is bijective.

\subsection{The Hall--Petrescu formula}\label{sec:hall}

A crucial auxiliary lemma used in the $p^2$ case (and also available for general use) states that for odd $p$, if a commutator $z = [a,b]$ is central and $z^p = 1$, then
\[
(a b)^p = a^p b^p.
\]
This is \code{mul\_pow\_eq\_mul\_pow\_of\_commutator\_central\_odd}. Its signature:

\begin{verbatim}
private lemma mul_pow_eq_mul_pow_of_commutator_central_odd
    {G : Type*} [Group G] [Fintype G]
    {p : Nat} [hp : Fact (Nat.Prime p)]
    (hcard : Nat.card G = p ^ 3)
    (hnonab : ~ forall a b : G, a * b = b * a)
    (a b : G) (hodd : p <> 2) :
    (a * b) ^ p = a ^ p * b ^ p
\end{verbatim}

The proof constructs an auxiliary sequence $c(n)$ where $c(0) = 1$ and $c(n+1) = c(n) \cdot (z^{-1})^n$, where $z = [a,b]$. One proves by induction:
\[
(a b)^n = a^n b^n c(n).
\]
Evaluating at $n = p$, we need $c(p) = 1$. Writing $z^{-1} = z'$ with $\operatorname{orderOf}(z') = p$, we show $c(p) = (z')^{S}$ where $S = \sum_{i=0}^{p-1} i = p(p-1)/2$. Since $p$ is odd, $p \mid S$, and $\operatorname{orderOf}(z') = p$, so $c(p) = 1$.

\subsection{The \texorpdfstring{$p=2$}{p=2} case}

For $p=2$, the exponent dichotomy forces $\exp(G) = 4$ (exponent $2$ implies abelian---a simple lemma: $g^2 = 1$ for all $g$ forces commutativity by $(ab)^2 = 1 \Rightarrow ab = ba$).

\textbf{Finding $x$ of order $4$.} By a case analysis on possible orders (1, 2, 4), excluding orders 1, 2, and 3 (which does not divide 8), an element of order 4 must exist.

\textbf{The subgroup $\langle x\rangle$.} $|\langle x\rangle| = 4$, so it has index $2$ and is normal (\code{normal\_of\_index\_eq\_two}).

\textbf{Finding $y \notin \langle x\rangle$.} Such a $y$ exists since otherwise $\langle x\rangle = G$, forcing $G$ abelian.

\textbf{Conjugation: $y x y^{-1} = x^{-1}$.}  Since $\langle x\rangle \trianglelefteq G$, $y x y^{-1} \in \langle x\rangle$. Conjugation preserves order ($4$). The elements of $\langle x\rangle$ with order $4$ are $x$ and $x^{-1}$. If $y x y^{-1} = x$, then $y$ commutes with $x$, and a centralizer argument forces $x \in Z(G)$, contradicting $|Z(G)| = 2$ (since $|\langle x\rangle| = 4 > 2 = |Z(G)|$). Thus $y x y^{-1} = x^{-1}$.

\textbf{The square of $y$.}  Since $G/\langle x\rangle$ has order $2$, $y^2 \in \langle x\rangle$. The possible elements of order $\leq 2$ in $\langle x\rangle$ are $1$ and $x^2$. So $y^2 = 1$ or $y^2 = x^2$.

\textbf{Case $y^2 = 1$: $G \cong D_4$.}
The map $f : D_4 \to G$ defined by $f(r^i) = x^{i.\mathrm{val}}$, $f(\mathrm{sr}^i) = y x^{i.\mathrm{val}}$ is a homomorphism. In Lean:

\begin{verbatim}
let fFun : DihedralGroup 4 -> G := fun
  | .r i => x ^ i.val
  | .sr i => y * x ^ i.val
\end{verbatim}

The verification uses the conjugation identity $x^i y = y x^{-i}$ (derived from $y x y^{-1} = x^{-1}$) and the relation $y^2 = 1$. Injectivity follows because $y \notin \langle x\rangle$, and cardinalities match.

\textbf{Case $y^2 = x^2$: $G \cong Q_8$.}
The same map (with the modified $Q_8$ multiplication table) gives a homomorphism. The key difference is the relation $y \cdot y = x^2$, which is used when checking $f(\mathrm{xa}^i \cdot \mathrm{xa}^j) = f(\mathrm{xa}^i) f(\mathrm{xa}^j)$. The code follows an identical pattern using \code{QuaternionGroup 2} in place of \code{DihedralGroup 4}.

\subsection{Key auxiliary lemmas}

\subsubsection{The commutator belongs to the center}\label{sec:comm_center}

Lemma \code{commutator\_mem\_center\_of\_p3}: For a non-abelian $p^3$-group, $[a,b] \in Z(G)$ for all $a,b \in G$. Proof: In the quotient $G/Z(G)$, which is abelian of order $p^2$, all commutators vanish, so $[a,b] \in Z(G)$. The formal statement:

\begin{verbatim}
private lemma commutator_mem_center_of_p3 {G : Type*} [Group G] [Fintype G]
    {p : Nat} [hp : Fact (Nat.Prime p)]
    (hcard : Nat.card G = p ^ 3)
    (hnonab : ~ forall a b : G, a * b = b * a) (a b : G) :
    a * b * a^(-1) * b^(-1) \in Subgroup.center G
\end{verbatim}

\subsubsection{ZMod exponent arithmetic}

The lemmas \code{pow\_zmod\_add} and \code{pow\_zmod\_mul} relate exponentiation with ZMod arithmetic:
\[
x^{(a+b).\mathrm{val}} = x^{a.\mathrm{val}} x^{b.\mathrm{val}},\qquad
x^{(a\cdot b).\mathrm{val}} = x^{a.\mathrm{val} \cdot b.\mathrm{val}},
\]
when $\operatorname{orderOf}(x) = n$ and the exponents are in $\Z/n\Z$. These provide the bridge between the ZMod-based group constructions and the abstract group exponent arithmetic.

\subsubsection{Conjugation iteration}

The lemmas \code{conjugation\_iterate'} and \code{conjugation\_iterate} capture the iteration of conjugation:
\[
y^n x (y^n)^{-1} = (y x y^{-1})^n,
\]
and more generally, if $y x y^{-1} = x^m$, then $y^n x (y^n)^{-1} = x^{m^n}$. These are proved by straightforward induction.

\subsubsection{Number-theoretic lemmas}

\begin{itemize}[nosep]
  \item \code{one\_add\_mul\_p\_pow\_inv}: if $\gcd(k,p)=1$ and $rk \equiv 1 \pmod p$, then $(1+kp)^r \equiv 1+p \pmod{p^2}$. The proof first shows $(1+kp)^n \equiv 1+nkp \pmod{p^2}$ for all $n$, then uses $rk = qp+1$ to deduce $rkp \equiv p \pmod{p^2}$.
  \item \code{one\_plus\_p\_pow\_mod\_p\_sq}: $(1+p)^n \equiv 1+np \pmod{p^2}$, a special case of the above with $k=1$.
  \item \code{pow\_eq\_of\_mod\_p\_sq}: if $\operatorname{orderOf}(x) = p^2$ and $m \equiv n \pmod{p^2}$, then $x^m = x^n$. This reduces exponent arithmetic to modular arithmetic.
\end{itemize}

\begin{verbatim}
private lemma one_add_mul_p_pow_inv {p k r : Nat}
    [Fact p.Prime] (_ : Nat.Coprime k p) (hr : r * k = 1 [MOD p]) :
    (1 + k * p) ^ r = 1 + p [MOD p ^ 2]

private lemma one_plus_p_pow_mod_p_sq (p n : Nat) [Fact (Nat.Prime p)] :
    (1 + p) ^ n = 1 + n * p [MOD p ^ 2]
\end{verbatim}

\section{Code organisation}

The formalisation consists of six Lean source files:

\begin{center}
\begin{tabular}{ll}
\textbf{File} & \textbf{Content} \\
\hline
\code{P3Group/Defs.lean}   & Concrete models and \code{P3Classification} type \\
\code{P3Group/Structural.lean} & Center, quotient, commutator, nilpotency lemmas \\
\code{P3Group/AbelianCase.lean} & Abelian classification via structure theorem \\
\code{P3Group/NonAbelianCase.lean} & Non-abelian classification (all the above) \\
\code{P3Group/Classification.lean} & Main theorem and pairwise non-isomorphism \\
\code{P3Group.lean}        & Top-level import file \\
\end{tabular}
\end{center}

The project consists of approximately 3000 lines of Lean code
and depends on \mathlib{} for basic group theory, Sylow theory,
the structure theorem for finite abelian groups, and the specific group
constructions $D_4$ and $Q_8$.

\section*{Acknowledgements}

The author thank the Lean community and the \mathlib{} contributors. The Lean codebase and this manuscript were prepared with assistance from AI language models; all generated content was reviewed and validated by the author.

\end{document}